\documentstyle[12pt,amssymb]{amsart}
\newcommand{\La}{\Lambda}

\newcommand{\GG}{{\cal G}}

\newcommand{\si}{\sigma}

\newcommand{\End}{\operatorname{End}}

\newcommand{\ed}{\qed\vspace{3mm}}

\newcommand{\lrar}[1]{\begin{picture}(50,10)(-25,-5)
\put(-25,0){\vector(1,0){50}}
\put(0,5){\makebox(0,0)[b]{\mbox{$#1$}}}
\end{picture}}

\newcommand{\ldar}[1]{\begin{picture}(10,50)(-5,-25)
\put(0,25){\vector(0,-1){50}}
\put(5,0){\mbox{$#1$}}
\end{picture}}

\newcommand{\ldrar}[1]{\begin{picture}(50,50)(-25,-25)
\put(-25,25){\vector(1,-1){50}}
\put(5,0){\mbox{$#1$}}
\end{picture}}

\newtheorem{thm}{Theorem}[section]
\newtheorem{prop}[thm]{Proposition}
\newtheorem{lem}[thm]{Lemma}

\theoremstyle{definition}
\newtheorem{defi}[thm]{Definition}
\newenvironment{rem}{\vspace{3mm}
\noindent {\bf Remark.}}{\vspace{3mm}}

\numberwithin{equation}{section}

\newcommand{\Pf}{\noindent {\it Proof}}
\newcommand{\Pic}{\operatorname{Pic}}
\newcommand{\Coh}{\operatorname{Coh}}
\renewcommand{\a}{\alpha}
\renewcommand{\b}{\beta}

\newcommand{\wt}{\widetilde}
\renewcommand{\mod}{\operatorname{mod}}

\renewcommand{\AA}{{\cal A}}

\newcommand{\FF}{{\cal F}}
\renewcommand{\GG}{{\cal G}}

\newcommand{\Z}{{\Bbb Z}}

\newcommand{\PP}{\cal P}

\newcommand{\Br}{\operatorname{Br}}
\newcommand{\GL}{\operatorname{GL}}
\newcommand{\PGL}{\operatorname{PGL}}

\newcommand{\ov}{\overline}

\newcommand{\ra}{\rightarrow}
\newcommand{\hra}{\hookrightarrow}

\newcommand{\id}{\operatorname{id}}
\newcommand{\G}{{\Bbb G}}

\newcommand{\ot}{\otimes}

\renewcommand{\O}{{\cal O}}

\newcommand{\D}{{\cal D}}
\newcommand{\De}{\Delta}

\newcommand{\PPic}{{\cal P}ic}
\newcommand{\sub}{\subset}

\title{Symplectic biextensions and a generalization
of the Fourier-Mukai transform}
\author{A. Polishchuk}

\begin{document}
\maketitle

Let $A$ be an abelian variety, $\hat{A}$ be the dual
abelian variety.
The Fourier-Mukai transform is an equivalence between
the derived categories of coherent sheaves $\D^b(A)$
and $\D^b(\hat{A})$. Notice that there is
a "symplectic" line bundle $L_A$ on $(\hat{A}\times A)^2$,
namely, $L_A=p_{14}^*\PP\ot p_{23}^*\PP^{-1}$ where $\PP$
is the Poincar\'e line bundle on $\hat{A}\times A$, such that
the standard embeddings $A\sub\hat{A}\times A$ and
$\hat{A}\sub\hat{A}\times A$ are "lagrangian" with
respect to $L_A$, i.e. the restrictions of $L_A$ to
$A^2$ and $\hat{A}^2$ are trivial (and they are maximal
with this property).
The purpose of this paper
is to establish an analogous equivalence of derived
categories for arbitrary
lagrangian subvarieties in an abelian variety $X$ equipped
with a line bundle $L$ over $X^2$ which satisfies some
properties similar to that of $L_A$ ($L$ should be
a {\it symplectic biextension}---see below).
Namely, with every lagrangian subvariety
$Y\sub X$ we associate a canonical element of the Brauer group
$e_Y\in \Br(X/Y)$ and consider the derived category
$\D^b(X/Y,e_Y)$ of
modules over the corresponding Azumaya algebra on $X/Y$.
It turns out that for every pair of lagrangian subvarieties
in $X$ there is an equivalence between these categories
generalizing the Fourier-Mukai transform.
The class $e_Y$ is trivial if and only if
the projection $X\ra X/Y$ splits, in this case
$\D^b(X/Y,e_Y)\simeq\D^b(X/Y)$.
This implies the "if" part of the
following   conjecture:   {\it  the  derived  categories  of
coherent sheaves on  abelian  varieties  $A$  and  $A'$  are
equivalent if and only if there is an isomorphism
$f:\hat{A}\times A\wt{\ra}\hat{A'}\times A'$ such that
$(f\times f)^*L_{A'}\simeq L_A$}.
In particular,
for any abelian variety $A$ and a symmetric homomorphism
$f:A\ra\hat{A}$ we construct an equivalence
$\D^b(A)\simeq\D^b(A/\ker(f_n))$ where $f_n=f|_{A_n}$
provided that $mn\ker(f)=0$ for some $m$ relatively prime
to $n$.

The construction
is   based   on   analogy   with  the  classical  theory  of
representations of the Heisenberg group of a symplectic vector space:
the categories $\D^b(X/Y,e_Y)$ are just different models of
the same "irreducible" representation of the {\it Heisenberg
groupoid}---a monoidal groupoid naturally attached to $(X,L)$.
The corresponding analogue of  Weil  representation
is studied in \cite{Weilrep}.

In what follows we consider varieties over an
algebraically closed field.

\section{Symplectic biextensions}

Let $X$ be an abelian variety.
A {\it biextension} of $X^2$ is a line bundle $L$ on $X^2$
together with isomorphisms
\begin{align*}
&L_{x+x',y}\simeq L_{x,y}\ot L_{x',y},\\
&L_{x,y+y'}\simeq L_{x,y}\ot L_{x,y'}
\end{align*}
--- this is a symbolic notation for isomorphisms
$(p_1+p_2,p_3)^*L\simeq p_{13}^*L\ot p_{23}^*L$ and
$(p_1,p_2+p_3)^*L\simeq p_{12}^*L\ot p_{13}^*L$ on $X^3$,
satisfying some natural cocycle conditions (see e.g. \cite{Breen}).

A  {\it  skew-symmetric  biextension}  of  $X^2$  is  a
biextension $L$ of $X^2$ together with
an   isomorphism   of   biextensions
$\phi:\si^*L\wt{\ra}  L^{-1}$, where   $\si:X^2\ra  X^2$  is  the
permutation of factors, and a trivialization
$\De^*L\simeq\O_X$ of $L$ over the diagonal
$\De:X\ra X^2$ compatible with $\phi$.

Every biextension $L$ of $X^2$ induces
a   homomorphism    $\psi_L:X\ra\hat{X}$    which    is
given on the level of points by $x\mapsto L_{x\times X}$.
If $L$ is skew-symmetric, then $\widehat{\psi_L}=-\psi_L$.
It is easy to see that $\psi_L=\psi_{L'}$ if and only if $L$ and $L'$
are isomorphic. Moreover, any skew-symmetric
homomorphism $\psi:X\ra\hat{X}$ defines a skew-symmetric
biextension  by  the formula $L(\psi)=(\psi\times\id)^*\PP$,
where $\PP$ is the  normalized  Poincar\'e  line  bundle  on
$\hat{X}\times X$, such that $\psi_{L(\psi)}=\psi$.
A skew-symmetric biextension $L$ is called {\it symplectic}
if $\psi_L$ is an isomorphism.

Let $Y\sub X$ be an abelian subvariety. Then $Y$  is  called
{\it isotropic}
with respect to $L$ if there is an isomorphism of skew-symmetric
biextensions  $L|_{Y\times Y}\simeq\O_{Y\times Y}$.
This is equivalent to the condition that the composition
$$Y\stackrel{i}{\ra} X\stackrel{\psi_L}{\ra}\hat{X}
\stackrel{\hat{i}}{\ra} \hat{Y}$$
is zero.
An isotropic subvariety $Y\sub X$ is called {\it lagrangian}
if the morphism $Y\ra \widehat{X/Y}$ induced by $\psi_L$ is an
isomorphism.
A  skew-symmetric  biextension  $L$  of  $X^2$  is called {\it
quasi-split} if there exists a lagrangian subvariety in $X$.
One can see easily that such a  biextension  is  necessarily
symplectic.
The simplest example of an abelian variety with a symplectic biextension
is $X=\hat{A}\times A$ for any abelian variety $A$ with
the biextension $L_A=P\ot\si^*P^{-1}$ where
$$P=p_{14}^*\PP\in\Pic(\hat{A}\times A\times \hat{A}\times A),$$
$\PP$ is the Poincar\'e line bundle on $A\times\hat{A}$.
A  symplectic  biextension  is  called  {\it split} if it is
isomorphic to this one.
Below we show how to construct all quasi-split symplectic biextensions.

Let $Y\sub X$ be a lagrangian subvariety.
Let us denote $A=X/Y\simeq\hat{Y}$ so that there is
an exact sequence
$$0\ra \hat{A}\ra X\stackrel{p}{\ra} A\ra 0,$$
such that $\psi_L|_{\hat{A}}=\hat{p}$.
The projection $p$ splits up to isogeny, that is there exists
a homomorphism $s:A\ra X$ such that $ps=n\id_A$.

\begin{lem}  One  can  always choose a section $s:A\ra X$ as
above such that \break $\hat{s}\psi_L s=0$.
\end{lem}

\Pf. Start with any $s$ as above and  then  replace  $n$  by
$2n^2$, and $s$ by $2ns-\hat{s}\psi_Ls$.
\ed

Choose $s:A\ra X$ as in lemma, and
let $\pi=(\id,s):\hat{A}\times A\ra X$ be the corresponding
isogeny. Then since $L|_{\hat{A}^2}$ and $(s\times  s)^*L$ are
trivial, it   is  easy  to  see  that  $\pi^*L\simeq
L_A^{\ot n}=p_{14}^*\PP^{\ot n}\ot p_{23}^*\PP^{\ot -n}$.
Thus,   every  quasi-split  symplectic  biextension descends
from the power of the split one. It remains to determine which
subgroups $\ker(\pi)\sub\hat{A}\times A$ can occur.

Let $P$ be any biextension of $X'\times X''$. Then the
restrictions of $P^{\ot n}$ on $X'_n\times X''$ and
$X'\times X''_n$ are canonically trivialized but these
trivializations
differ over $X'_n\times X''_n$ by a bilinear morphism
$e_n(P):X'_n\times X''_n\ra\G_m$. In the case of the Poincar\'e
line bundle $\PP$ over $\hat{A}\times A$ this construction
gives   a  canonical  perfect  pairing  $e_n:\hat{A}_n\times
A_n\ra\G_m$.
In our situation the canonical trivializations of  $L_A^{\ot
n}$ over $(\hat{A}_n\times A_n)\times (\hat{A}\times A)$ and
$(\hat{A}\times A)\times (\hat{A}_n\times A_n)$ differ
over $(\hat{A}_n\times A_n)^2$ by a bilinear morphism
$e_n(L_A):(\hat{A}_n\times A_n)^2\ra\G_m$ which is given by the formula
\begin{equation}
e_n(L_A)((\xi,x),(\xi',x'))=e_n(\xi,x')e_n(\xi',x)
\end{equation}
where $x,x'\in A_n$, $\xi,\xi'\in\hat{A}_n$.

By definition $\ker{\pi}$ is the graph of a morphism
$\phi:A_n\ra  \hat{A}_n$  induced  by  $s$.
Now the biextension $L_A^{\ot n}$ descends to $X$ if and only
if   there  exist  trivializations (as a biextension) of
$L_A^{\ot  n}$  over
$\ker{\pi}\times  (\hat{A}\times  A)$  and   $(\hat{A}\times
A)\times\ker{\pi}$ which are compatible over $(\ker{\pi})^2$.
Since  such  trivializations  are  unique  they
coincide   with   the   restrictions   of   the    canonical
trivializations above. Hence, the descent condition is that
$\ker(\pi)$ is isotropic with respect to $e_n(L_A)$ which
means  that  $\phi:A_n\ra\hat{A}_n$  is  skew-symmetric with
respect to $e_n$, that is $\widehat{\phi}=-\phi$.

Thus, any quasi-split symplectic biextension arises from
a   pair   $(A,\phi)$,   where   $\phi:A_n\ra\hat{A}_n$  is  a
skew-symmetric morphism, as described above.
It is easy to see that  if  we
change $\phi$ by $\phi+f_n$ where $f_n$ is the restriction
of  a  symmetric  homomorphism  $f:A\ra\hat{A}$ to $A_n$ (then
$f_n$ is automatically skew-symmetric), then we get isomorphic
symplectic biextensions---this corresponds to  a  change  of
an isotropic morphism  $s:A\ra X$. Also, one can change $n$ by
$nm$ and $\phi$ by the composition
$$A_{nm}\stackrel{m}{\ra}A_n\stackrel{\phi}{\ra}\hat{A}_n\ra
\hat{A}_{nm},$$
so that the corresponding symplectic biextension will be the
same. However, this doesn't exhaust examples of pairs
$(A,\phi)$ giving isomorphic biextensions.
For example, it is easy to see that $A/\ker(\phi)=s(A)\sub X$
is a lagrangian subvariety in $X$, $X/s(A)\simeq
\hat{A}/\phi(A_n)$, and the biextension associated with the
pair    $(A,\phi)$    corresponds    also    to  the   pair
$(\hat{A}/\phi(A_n),\psi)$ where $\psi$ is the composition
$$\psi:(\hat{A}/\phi(A_n))_n\ra
\phi(A_n)\stackrel{\phi^{-1}}{\ra} A_n/\ker(\phi)\ra
(A/\ker(\phi)_n).$$

These considerations lead to the following theorem.

\begin{thm}\label{pairs}
Let $L$ be a symplectic biextension of $X^2$.
For any lagrangian subvariety $Y\sub X$
there exists a lagrangian subvariety $Z\sub X$ such that
$Y\cap Z$ is finite. Any pair of lagrangian subvarieties
$(Y,Z)$ in $X$ such that $Y\cap Z$ is finite, is isomorphic
to the pair $(\hat{A},A/\ker(\phi))$ in
$\hat{A}\times   A/(\phi,\id)(A_n)$   with   its   canonical
symplectic biextension for some abelian variety $A$ and
a skew-symmetric  homomorphism  $\phi:A_n\ra\hat{A}_n$.
\end{thm}

\Pf . The first assertion is clear. To prove the second
we should start with the lagrangian subspace $Y\sub X$
in the above argument and choose a splitting of $p:X\ra X/Y$
up to isogeny which factors through $Z$. More precisely,
let $f:Z\ra X/Y$ be the restriction of $p$ to $Z$. Choose an
isogeny $g:X/Y\ra Z$ such that $fg=n\id_{X/Y}$. Then
the composition of $g$ with the embedding of $Z$ in $X$ gives
a lagrangian morphism $s:X/Y\ra X$ such that $ps=n\id_{X/Y}$.
Now we get an isogeny $\hat{A}\times A\ra X$ as above (where
$A=X/Y$) such that $Y$ and $Z$ are the images of  $\hat{A}$  and  $A$
respectively, which finishes the proof.
\ed

Let us give an example of a quasi-split symplectic
biextension which is not split.
Let $A$ be a principally polarized abelian variety with
$\End(A)=\Z$. Then there is a symplectic isomorphism
$\phi_n:A_n\ra\hat{A}_n$ such that for every symmetric morphism
$f:A\ra  \hat{A}$  the corresponding morphism $f|_{A_n}$ is proportional to
$\phi_n$. Now if $\dim(A)>1$  we  can  choose  a  symplectic
morphism $\phi:A_n\ra\hat{A}_n$
which  is  not  proportional  to $\phi_n$. It is easy to see
that the corresponding symplectic biextension of $X^2$, where
$X=\hat{A}\times A/(\phi,\id)(A_n)$, is not split.

\section{Representations of the Heisenberg groupoid}

Let   $X$  be  an  abelian  variety,  $L$  be  a  symplectic
biextension of $X^2$. Throughout this section we
assume that there exists a biextension
$P$ of $X^2$ such that $L\simeq  P\ot\si^*P^{-1}$
(an isomorphism of skew-symmetric biextensions).  This  is
equivalent to the condition $\psi_L=f-\hat{f}$ for some
$f:X\ra  \hat{X}$.  For  example,
the quasi-split biextension associated with a
pair $(A,\phi)$, where $\phi:A_n\ra\hat{A}_n$ and $n$ is odd,
satisfies this condition.

\begin{defi}         The         Heisenberg         groupoid
$H(X)=H(X,P)$ is the stack of monoidal groupoids
such that $H(X)(S)$ for a scheme $S$ over $k$
is the monoidal groupoid generated by
the central subgroupoid $\PPic(S)$ of $\G_m$-torsors on $S$ and the symbols
$T_x$, $x\in X(S)$ with the composition law
$$T_x\circ T_{x'}= P_{x,x'} T_{x+x'}.$$
In other words, objects of $H(X)(S)$ are pairs $(M,x)$ where
$M$ is a line bundle over $S$, $x\in X(S)$.
A morphism $(M,x)\ra(M',x')$ exists only if $x=x'$ and is
given by an isomorphism $M\ra M'$. The composition law
is defined by the formula
$$(M,x)\circ (M',x')=(P_{x,x'}\ot M\ot M', x+x').$$
Denoting $T_x=(\O_S,x)$ we recover the above relation.
\end{defi}

If we replace $P$ by $P'=P\ot\La(M)$ for some line bundle
$M$ on $X$ trivialized along the zero section,
where $\La(M)=(p_1+p_2)^*M\ot p_1^*M^{-1}\ot p_2^*M^{-1}$
(see e.g. \cite{Breen})
we get an equivalent Heisenberg groupoid. The equivalence
$H(X,P)\ra H(X,P')$ is defined by the functor which is the identity on
$\PPic(S)$ and sends $T_x$ to $M_x^{-1}T_x$.
Since any symmetric biextension of $X^2$ has form $\La(M)$
this shows that up to a non-unique equivalence the Heisenberg
groupoid doesn't depend on a choice of $P$ such that $L=P\ot\si^*P$.

\begin{rem} One can see easily that the Heisenberg groupoid
can be considered as an
extension of the group scheme $X$ by the stack of line bundles
in the sense of Deligne (see \cite{Des}), namely, we associate to each
point $x\in X(S)$ the trivial gerb of line bundles,
and the composition is given by the formula above.
\end{rem}

The  Heisenberg  groupoid $H(\hat{A}\times A)$ corresponding
to a split biextension
is generated by the Picard subgroupoid $\PPic$ and symbols $T_x$, $T_y$ where
$x\in \hat{A}$, $y\in A$ with the following defining relations:
\begin{align*}
&T_xT_{x'}=T_{x+x'},\\
&T_yT_{y'}=T_{y+y'},\\
&T_y T_x=\langle  x,y\rangle T_x T_y.
\end{align*}

Is is easy to see (see e.g. \cite{Weilrep}) that
the map $T_y\mapsto t_y^*$, $T_x\mapsto \cdot\ot\PP_x$
defines an action of $H(\hat{A}\times A)$ on $\D^b(A)$,
where $t_y:A\ra A$ is the translation by $y\in A$,
$\PP_x=\PP|_{x\times A}$ for $x\in\hat{A}$.
Below we construct an analogous action for an
arbitrary isotropic subvariety of an abelian variety with
a symplectic   biextension.

Let $Y\sub X$ be an isotropic subvariety. Then
$P|_{Y\times Y}$ has a natural structure of a symmetric biextension.

\begin{defi}  A  pair  $(Y,\a)$,  where  $Y$  is  an isotropic abelian
subscheme of $X$ with respect to $L$ and $\a$ is a line bundle
on $Y$ with fixed trivialization along the zero section, is  called
isotropic if
an isomorphism of symmetric biextensions of $Y\times Y$ is given:
$$\La(\a)\simeq P|_{Y\times Y}$$
which can be written symbolically as
$P_{y,y'}=\a_{y+y'}\a_y^{-1}\a_{y'}^{-1}$ for
$y,y'\in Y$.
\end{defi}

For  any  isotropic  subvariety  $Y\sub X$ there exists $\a$
such that the pair $(Y,\a)$ is isotropic.

\begin{defi}  For  an  isotropic pair $(Y,\a)$ we define
$\FF(Y,\a)$ as the category of pairs $(\AA,a)$ where
$\AA\in\D^b(X)$, $a$ is an isomorphism in $\D^b(Y\times X)$:
\begin{equation}\label{Schrsp}
a:(ip_1+p_2)^*\AA\wt{\ra} P^{-1}|_{Y\times X}
\ot p_1^*\a^{-1}\ot p_2^*\AA
\end{equation}
where $i:Y\hra X$ is the embedding, such that $(e\times\id)^*a=\id$.
This isomorphism can be written
symbolically as follows:
$$a_{y,x}:\AA_{y+x}\wt{\ra}P^{-1}_{y,x}\a^{-1}_y\AA_x$$
where $y\in Y$, $x\in X$.
These data should satisfy
the following cocycle condition:
$$a_{y_1+y_2,x}=a_{y_2,x}\circ a_{y_1,y_2+x}:\AA_{y_1+y_2+x}
\ra P^{-1}_{y_1,x+y_2} P^{-1}_{y_2,x}\a^{-1}_{y_1}\a^{-1}_{y_2}\AA_x
\simeq P^{-1}_{y_1+y_2,x}\a^{-1}_{y_1+y_2}\AA_x,$$
or in standard notation
$$(p_1+p_2,p_3)^*a=(p_2,p_3)^*a\circ (p_1,ip_2+p_3)^*a$$
in $\D^b(Y\times Y\times X)$.
The morphisms between such pairs are morphisms
between the corresponding objects in $\D^b(X)$ commuting with the
isomorphisms in (\ref{Schrsp}).
\end{defi}

It is easy to see that the category
$\FF(Y,\a)$ is equivalent to $\D^b(X/Y)$ provided the
projection $p:X\ra X/Y$ has a section $s:X/Y\ra X$.
However, in general this is not true: one encounters
some twisted versions of $\D^b(X/Y)$ considered in the next
section.

There is a natural action of the Heisenberg groupoid $H(X)$
on the category $\FF(Y,\a)$ such that an object $(M,x)$ acts
by the functor
$$\AA\mapsto M\ot P|_{X\times x}\ot t_x^*(\AA).$$
In the case $X=\hat{A}\times A$, $Y=\hat{A}\sub X$
this action coincides with the action of $H(\hat{A}\times A)$ on
$\D^b(A)\simeq\FF(\hat{A})$ mentioned above.

By analogy with the classical Heisenberg group it is natural to ask
when these representations are irreducible in some sense.
More precisely, for the construction of Weil representation it is relevant
to know that all intertwining operators from Schr\"odinger representation
to itself are proportional to the identity. As shown in \cite{Weilrep}
certain analogue of this property holds for the action of $H(\hat{A}\times A)$
on $\D^b(A)$. One can treat the case of an arbitrary lagrangian
subvariety similarly, however, we don't need this result.

\section{Modules over Azumaya algebras}

We begin this section by briefly recalling the various ways
to speak about the category of coherent modules over a scheme
$S$  "twisted"  by  an  element $e\in H^2(S,\G_m)$: via Cech
cocycles, gerbs, and Azumaya algebras. The  simplest  way  to
define such a category is to fix an open covering $(U_i)$ of
$S$ (say, in flat topology)
such that $e$ is represented by a Cech cocycle
$\a_{ijk}\in\O^*(U_{ijk})$ where
$U_{ijk}=U_i\times_S U_j\times_S U_k$.
Then we define $\Coh(S,\a)$ as the category of collections
$(\FF_i)$ of coherent sheaves on $U_i$ together with a
system of isomorphisms $f_{ij}:\FF_i\ra\FF_j$ over
$U_{ij}=U_i\times_S  U_j$ (such that $f_{ji}=f_{ij}^{-1}$)
satisfying  the  twisted  cocycle
condition: $f_{jk}f_{ij}=\a_{ijk}f_{ik}$ over $U_{ijk}$.
It is easy to see that up to equivalence this category
depends only on the cohomology class of $\a$.
The more abstract way to define this category (which doesn't
involve  a choice of covering) is to represent $e$ by a
$\G_m$-{\it gerb}. Recall that a $\G_m$-gerb is a stack of
groupoids  $\GG$  such  that  locally  there  is  a   unique
isomorphism  class of objects of $\GG$ and the automorphism
group  of  any  object  is  $\G_m$.  Equivalence  classes  of
$\G_m$-gerbs over $S$ are in bijective correspondence with
$H^2(S,\G_m)$. Now consider the category of
representations of $\GG$, i.e. the category of functors of
stacks $\GG\ra\Coh(S)$ where $\Coh(S)$ is the stack of
coherent sheaves. Choosing an open covering and a collection
of objects $V_i\in\GG(U_i)$ we arrive to the Cech description
above. Sometimes $e$ is represented
by a sheaf of Azumaya algebras $\AA$ over $S$. Then locally
$\AA$  is  isomorphic to a matrix algebra of rank $n^2$ over
$S$. Now let $\GG(\AA)$ be the $\G_m$-gerb  of  representations  of
$\AA$ in locally free $\O_S$-modules of rank $n$. Then it is
easy  to  see that $\GG(\AA)$ represents the same cohomology
class $e\in H^2(S,\G_m)$ and the categories of representations
of $\GG(\AA)$ and $\AA$  in  coherent  sheaves  on  $S$  are
equivalent. By abuse of notation we denote all these equivalent
categories by $\Coh(S,e)$.

Let  $E\ra  S$  be  a  $K$-torsor where $K$ is a finite flat
commutative group
scheme over $S$, let $0\ra\G_m\ra G\ra K\ra 0$ be a  central
extension of $K$. Then it defines an element $e(G,E)\in H^2(S,e)$
such that the category of $G$-equivariant  coherent  sheaves
on $E$ of weight 1 is equivalent to $\Coh(S,e(G,E))$. Here
a weight of a $G$-equivariant coherent sheaf is defined as
the weight of the induced $\G_m$-equivariant sheaf.
Indeed,  consider  the  gerb  $\GG(G,E)$  of  liftings  of $E$ to
$G$-torsors (an object of $\GG(G,E)$ over $U\ra S$ is a $G$-torsor
$\wt{E}$ over $U$ together with an isomorphism of $K$-torsors
$\wt{E}/\G_m\simeq E$). Then we claim that  the
category of weight-1 $G$-equivariant sheaves is equivalent to
the category of representations of $\GG(G,E)^{op}$ which is
$\Coh(S,e)$ where $e$ is the inverse of the cohomology class of $\GG(G,E)$.
To see this note that a lifting of $E$ to a  $G$-torsor  can
be considered as a weight-1 $G$-equivariant line bundle $L$ over $E$.
A  choice   of  such  bundle over $E_U$ defines  the equivalence
$\FF\mapsto \FF\ot  L^{-1}$ of the category of weight-1  $G$-equivariant
sheaves with the category of $K$-equivariant sheaves on $E_U$,
hence with $\Coh(U)$ depending contravariantly on $L$,
whereas the assertion.
The class $e(G,E)$ is trivial if and only if there is a global
object  of  $\GG(G,E)$,  i.e.  a  global lifting of $E$ to a
$G$-torsor. Also it is easy to see that $e(G,E)$ depends
biadditively  on the pair of classes $[G]\in H^2(K,\G_m)$,
$[E]\in H^1(S,K)$.

We'll apply this in the particular case when $S=A$
is an abelian variety, $p:E\ra A$ is  an  isogeny  of  abelian
varieties, so that $E$ can be considered as a $K$-torsor where
$K=\ker(p)$.  Then  for  any central extension $\pi:G\ra K$ by
$\G_m$ the previous construction gives  a  class  $e(G,E)\in
H^2(A,\G_m)$ which is an obstruction for existence of a line
bundle  $M$  over  $E$ such that $K\sub K(M)$ and $G$ is the
restriction of Mumford's extension $G(M)\ra K(M)$ to $K$
(see \cite{Mum}).
Let $\rho:G\ra\GL_n$ be a weight-1 representation of $G$,
$\ov{\rho}:K\ra\PGL_n$ be the corresponding projective
representation of $K$. Then the $PGL_n$-torsor $E_{\ov{\rho}}$
on $A$ obtained as the push-forward of $E$ by $\ov{\rho}$
gives rise to an Azumaya algebra with the class $e(G,E)$.
Consider $G$ as a $\G_m$-torsor over $K$ so that
$G_u=\pi^{-1}(u)$ for $u\in K$. Let us denote by $\O_K(G)$
the corresponding line bundle over $K$.
Then a weight-1  $G$-equivariant  sheaf
on $E$ can be described by the following data:
a coherent sheaf $\FF$ on $E$ and an isomorphism
over $K\times E$:
\begin{equation}\label{equiv}
p_1^*\O_K(G)\ot p_2^*\FF\wt{\ra} (ip_1+p_2)^*\FF
\end{equation}
where $i:K\ra E$ is the inclusion, satisfying the natural cocycle
condition. The above construction gives an equivalence of
this category with $\Coh(A,e(G,E))$.

We need also a derived category version of this equivalence.
The slight difficulty is that derived categories of coherent
sheaves don't glue well in any of standard topologies.
However, as shown in the Appendix,
the descent formalism for
finite flat morphisms extends to derived categories.
This allows to rephrase the definition of $\Coh(S,e)$
(e.g. in Cech version) for a class $e$ which is killed by
a finite flat morphism $S'\ra S$ into a description of
the corresponding derived category $\D^b(S,e)$. Similarly,
one can describe the derived category of weight-1 $G$-equivariant
sheaves above as the category of objects $\FF\in\D^b(E)$ with
isomorphisms (\ref{equiv}) satisfying the cocycle condition
and to show that it is equivalent to $\D^b(A,e(G,E))$.

Let   $X$  be  an  abelian  variety,  $L=P\ot\si^*P^{-1}$
be  a  symplectic
biextension of $X^2$, $(Y,\a)$ be an isotropic pair.

\begin{prop}
There is a canonical class $e(Y)\in
H^2(X/Y,\G_m)$  such  that  the category $\FF(Y,\a)$ defined
in the previous section is equivalent to $\D^b(X/Y,e(Y))$.
\end{prop}

\Pf . Choose a homomorphism of abelian varieties  $s:Z\ra  X$  and  a  line
bundle $\b$ on $Z$ such  that  the
restriction  of  the  composition $ps:Z\ra X/Y$ is an
isogeny and there is an isomorphism of biextensions
of $s^{-1}(Y)\times Z$
\begin{equation}\label{inters}
(s\times s)^*P|_{s^{-1}(Y)\times Z}\simeq
\La(\b)|_{s^{-1}(Y)\times Z}.
\end{equation}
For  example,  one can take $Z=X/Y$ and $s':Z\ra X$ such that
$ps'=n\id_{X/Y}$, then $s^{\prime -1}(Y)=(X/Y)_n$ and
$(ns'\times ns')^*P|_{(X/Y)_{n^2}\times X/Y}$ is a trivial biextension
(see \cite{Breen}, 4.2), hence we can take $s=ns'$ and $\b=\O_Z$.

Then $\La((s|_{s^{-1}(Y)})^*\a\ot \b^{-1}|_{s^{-1}(Y)})$ is a trivial
biextension, hence the $\G_m$-torsor
$\b|_{s^{-1}(Y)}\ot (s|_{s^{-1}(Y)})^*\a^{-1}$
defines a central extension $G$ of $s^{-1}(Y)$ by $\G_m$.
It is easy to see that
the class $e(G,Z)\in H^2(X/Y,\G_m)$ defined
above doesn't depend on a choice of $\a$ such that the
pair $(Y,\a)$ is isotropic. We claim also that it
doesn't depend on a choice of $Z$ and $\b$.
Indeed, if we change $\b$ by $\a'_Z=\b\ot\gamma$ where
$\La(\gamma)|_{s^{-1}(Y)\times Z}$ is trivial, then the new
central extension is the sum of $G$ and the restriction of
Mumford's extension $G(\gamma)\ra K(\gamma)$ to $s^{-1}(Y)$.
But  $\gamma$  has  a  natural
structure of $G(\gamma)$-equivariant line bundle, hence
$e(G(\gamma),Z)=0$.
Also, it is easy to see that $e(G,Z)$ doesn't change if we
replace $s:Z\ra X$ by the composition of $s$ with an isogeny
$Z'\ra  Z$.  It  remains to check that $e(G,Z)$ is invariant
under the change $s'=s+f$ where  $f:Z\ra  Y$  is  any
homomorphism. In this case $s^{\prime -1}(Y)=s^{-1}(Y)$ and
$$(s'\times s')^*P|_{s^{-1}(Y)\times Z}\simeq
\La(\a'_Z)|_{s^{-1}(Y)\times Z}$$
where $\a'_Z=\b\ot f^*\a\ot (f,s)^*P$. On the other hand
$$(s'|_{s^{-1}(Y)})^*\a\simeq (s|_{s^{-1}(Y)})^*\a\ot
(f^*\a\ot (f,s)^*P)|_{s^{-1}(Y)}$$
so  that  the corresponding central extension of $s^{-1}(Y)$
is the same.

Given   an   object   $(\AA,a)$   of   $\FF(Y,\a)$   where
$\AA\in\D^b(X)$, $a$ is an isomorphism (\ref{Schrsp}),
we  can  consider  $s^*\AA\in\D^b(Z)$.  Then  $a$ induces an
isomorphism
$$a_{s(u),s(z)}:s^*\AA_{u+z}\wt{\ra}(s\times s)^*P^{-1}_{u,z}
\a^{-1}_{s(u)} s^*\AA_z$$
where $u\in s^{-1}(Y)$, $z\in Z$, satisfying the usual cocycle
condition. Let $F(\AA)=s^*\AA\ot\b$, then $a_{s(u),s(z)}$
together with (\ref{inters}) gives an isomorphism
\begin{equation}
F(a):F(\AA)_{u+z}\wt{\ra}\b_u\ot \a^{-1}_{s(u)}\ot F(\AA)_z
\end{equation}
where $u\in s^{-1}(Y)$, $z\in Z$. In other words,
$(F(\AA),F(a))$ can be considered as a weight-1 $G$-equivariant object
of $\D^b(Z)$. This gives the required equivalence as one can
check  applying  Theorem  A  of  Appendix  to  the  morphism
$Y\times Z\ra X$ and the trivial descent for the projection
$Y\times Z\ra Z$.
\ed

\begin{rem} If $(Y,\a)$ and $(Z,\b)$ are isotropic pairs the
$\G_m$-torsor $\b|_{Y\cap Z}\ot\a|^{-1}_{Y\cap Z}$
defines a central extension $G$ of $Y\cap Z$ by $\G_m$,
which  in case when $Y$ and $Z$ are lagrangian and $Y\cap Z$
is finite gives the class $e_Y\in H^2(X/Y,\G_m)$.
The corresponding commutator form on $Y\cap Z$
measures the difference between the symmetric structures
on $P|_{(Y\cap Z)^2}$ restricted from $Y^2$ and $Z^2$.
In other words, this is the standard symplectic form
associated with the biextension $L|_{Y\times Z}$
measuring the difference between two trivializations of
$L|_{(Y\cap Z)^2}$ restricted from $Y\times (Y\cap Z)$
and $(Y\cap Z)\times Z$.
\end{rem}

By  definition  the  class $e(Y)$ vanishes if the projection
$X\ra X/Y$ splits. It turns out that if $Y$ is lagrangian then
the converse is also true.

\begin{prop} Let $Y\sub X$ be a lagrangian subvariety.
If $e_Y=0$ then the projection $X\ra X/Y$ splits.
\end{prop}

\Pf . According to Theorem \ref{pairs} we can assume
that $X=\hat{A}\times A/(\phi,\id)(A_n)$,
$Y=A/\ker(\phi)\sub X$ for an
abelian variety $A$ and a skew-symmetric homomorphism
$\phi:A_n\ra\hat{A}_n$. Now we can take $Z=\hat{A}\sub X$
in  the definition of $e_Y$. The kernel of the projection
$Z\ra X/Y$ is $\phi(A_n)\sub\hat{A}$ and the commutator form
of its central extension considered above is (up to sign)
$$e(\phi(x),\phi(y))=e_n(\phi(x),y)$$
where $x,y\in A_n$.
The triviality of $e_Y$ implies
that there exists a symmetric homomorphism $g:\hat{A}\ra A$
such that $\phi(A_n)\sub\ker(g)$ and $e=e^g|_{\phi(A_n)^2}$
where $e^g$ is the standard symplectic form on $\ker(g)$.
In other words, the following equality holds:
$$e_n(\phi(x),y)=e_n(\phi(x),g(n^{-1}\phi(y)))$$
for all $x,y\in A_n$, which implies that
$y-g(n^{-1}\phi(y))\in\ker(\phi)$ for $y\in A_n$.
Note  that  $x\mapsto   g(n^{-1}x)\mod(\ker(\phi))$   is   a
well-defined homomorphism $\hat{A}\ra A/\ker(\phi)$ since
$g(A_n)\sub    \ker(\phi)$    (which    is   obtained   from
$\phi(A_n)\sub\ker(g)$ by duality). Thus, a homomorphism
$$\hat{A}\times        A\ra         A/\ker(\phi):(x,y)\mapsto
y-g(n^{-1}x)\mod(\ker(\phi))$$
descends to a homomorphism $X=\hat{A}\times A/(\phi,\id)(A_n)\ra
A/\ker(\phi)=Y$ splitting the embedding $Y\ra X$.
\ed

\section{Intertwining functors}

Let  $X$ be an abelian variety with a symplectic biextension
$L$ of $X^2$.
In this section we construct an equivalence
of $H(X)$-representations
$\FF(Y,\a)\simeq\FF(Z,\b)$
for isotropic pairs $(Y,\a)$ and
$(Z,\b)$ such that $Y$ and $Z$ are lagrangian.

The idea is to mimic the classical construction.
Namely, consider the functor of ``integration over $Z$''
$$R:\FF(Y,\a)\ra\FF(Z,\b):
\AA\mapsto p_{2*}(P|_{Z\times X}\ot p_1^*\b\ot (ip_1+p_2)^*\AA).$$
The following symbolic notation stresses the analogy with the classical case:
$$R(\AA)_x=\int_{Z} P_{z,x}\b_z\AA_{z+x}dz.$$
It easy to check that $R(\AA)$ has a natural structure of
an object of $\FF(Z,\b)$:
\begin{align*}
&R(\AA)_{z'+x}=\int_{Z} P_{z,z'+x}\b_z\AA_{z+z'+x}dz\simeq
\int_{Z} P_{z,x}\b_{z+z'}\b_{z'}^{-1}\AA_{z+z'+x}dz\simeq\\
&\int_{Z} P_{z-z',x}\b_z\b_{z'}^{-1}\AA_{z+x}dz\simeq
P^{-1}_{z',x}\b_{z'}^{-1}R(\AA)_x
\end{align*}
--- here we used the isomorphism $P|_{Z^2}\simeq\La(\b)$ and
the change of variable $z\mapsto z-z'$.
It is also clear that $R$
commutes with the action of $H(X)$.
In the classical theory in order to get an invertible intertwining operator
one should replace the integration over $Z$ by the integration over
$Z/Y\cap Z$ in the above formula. This doesn't work literally in our context---
it turns out that in the correct definition
one eliminates the "excess" integration over the connected component
of $Y\cap Z$, and over a "largangian half" of the group of connected
components of $Y\cap Z$. Instead of working out the case when
$\dim(Y\cap Z)>0$ we use the following simple lemma which
allows to avoid it.

\begin{lem}  For  any  pair  $Y$  and  $Z$   of   lagrangian
subvarieties  of  $X$  there  exists a lagrangian subvariety
$T\sub X$ such that the intersections $Y\cap T$ and $Z\cap T$
are finite.
\end{lem}

\Pf . We can work in the category of  abelian  varieties  up  to
isogeny.  We have an isogeny $X\sim Y\times \hat{Y}$ and
$Z/Y\cap Z\sub Y\times \hat{Y}$
is isogenic to the graph of a symmetric morphism
$$g: Z/Y\cap Z\ra Y/Y\cap Z\sim \widehat{Z/Y\cap Z}.$$
Let $K\sim\ker(g)$.
We have a decomposition $Z/Y\cap Z\sim K\times K'$
such that $g$ is given by a symmetric isogeny $K'\ra\hat{K'}$.
Now   let   $\hat{Y}\sim  Z/Y\cap  Z\times  K''\sim  K\times
K'\times K''$. Let us define a symmetric morphism
$f:\hat{Y}\ra Y$ to be a symmetric isogeny on $K$ and zero
on two other factors. Then we can take the graph of  $f$  to
be $T$.
\ed

Thus, we may assume that $Y\cap Z$ is finite.
We have a natural central extension $G$ of $Y\cap Z$ by $\G_m$
given by the $\G_m$-torsor $\b|_{Y\cap Z}\ot\a^{-1}|_{Y\cap Z}$
such that $\FF(Y,\a)$ is
equivalent  to  the  category  of  weight-1  $G$-equivariant
objects of $\D^b(Z)$, while $\FF(Z,\b)$---to that of weight-1
$G^{-1}$-equivariant objects of $\D^b(Y)$, where $G^{-1}$ is
the   inverse   central  extension of $Y\cap Z$ (given  by  the  inverse
$\G_m$-torsor).
Let $e$ be the commutator form of $G$. Choose  a  lagrangian
subgroup  $H\sub  Y\cap Z$ with respect to $e$ and a
trivialization of the central extension $G$ over $H$  (which
is the same as a lifting of $H$ to a subgroup in $G$).
Then we can define the reduced functor
$$\ov{R}:\FF(Y,\a)\ra\FF(Z,\b):
\ov{R}(\AA)_x=\int_{Z/H} P_{z,x}\b_z\AA_{z+x}dz.$$
To give a meaning to this notice that
an object $P_{z,x}\b_z\AA_{z+x}\in\D^b(Z\times X)$  descends
canonically  to  an  object  of  $\D^b(Z/H\times  X)$ (use the
additional data on $\AA\in\FF(Y,\a)$ and the
isomorphism $\a|_H\simeq\b|_H$). As  above  it  is  easy  to
check  that  $\ov{R}(\AA)$  has  a natural structure of an object of
$\FF(Z,\b)$ and $\ov{R}$ commutes with $H(X)$-action.

\begin{thm}\label{main}
The  functor  $\ov{R}$  is  an  equivalence  of
categories.
\end{thm}

\Pf . First let us rewrite $\ov{R}$ as the functor from
the category of weight-1 $G$-equivariant objects of $\D^b(Z)$
to that of weight-1 $G^{-1}$-equivariant objects of $\D^b(Y)$
using   the   equivalences  defined  above. Recall that an
equivalence of the first category with $\FF(Y,\a)$ is given
by a functor $F_Y$ which associates to $\AA\in\FF(Y,\a)$
the $G$-equivariant object $\AA|_Y\ot\b\in\D^b(Z)$, while
the second equivalence is induced by
$F_Z:\FF(Z,\b)\ra\D^b(Y):\AA'\mapsto\AA'|_Y\ot\a$.
Now for $\AA\in\FF(Y,\a)$ we have
\begin{align*}
F_Z(\ov{R}(\AA))_y &=
\a_y\int_{Z/H} P_{z,y}\b_z\AA_{z+y}dz\simeq
\a_y\int_{Z/H}  P_{z,y}\b_z P_{y,z}^{-1}\a_y^{-1}\AA_z
dz\simeq\\
&\simeq\int_{Z/H} L_{z,y} F_Y(\AA)_z dz
\end{align*}
The latter integral should be understood in the same sense
as above: the $G$-equivariance data on $F_Y(\AA)$ allow
to    descend   $L_{z,y} F_Y(\AA)_z$   to   an   object   of
$\D^b(Z/H\times Y)$. Notice  that $G$-equivariance
data on an object $\GG\in\D^b(Z)$ includes the descent data
for the projection $Z\ra Z/H$, so that $G$-equivariant
objects  of  $\D^b(Z)$  can  be  considered  as  objects  of
$\D^b(Z/H)$   with   some    additional    data.
More precisely, the isomorphism
$$\GG_{z+u}\simeq G_u L_{z,u}\GG_z$$
where $u\in Y\cap Z$, induced by the $G$-equivariance data and
the trivialization of $L_{z,u}$ commutes with the
descent data for $Z\ra Z/H$, so it induces an isomorphism
of descended objects on $Z/H\times (Y\cap Z)$
$$\ov{\GG}_{z+u}\simeq G_u \ov{L}_{z,u}\ov{\GG}_z$$
--- these are the additional data for $\ov{\GG}\in\D^b(Z)$.
Similar,
we can consider $G^{-1}$-equivariant objects of $\D^b(Y)$
as   objects   of  $\D^b(Y/H)$  with  additional  data. It is
easy to see that biextension $L_{z,y}$ of $Z\times Y$ descends to a
biextension $\ov{L}$ of $Z/H\times Y/H$ which induces an isomorphism
$Z/H\wt{\ra}\widehat{Y/H}$. Thus, $\ov{R}$ is compatible with
the  Fourier-Mukai transform  $\D^b(Z/H)\ra\D^b(Y/H)$ via
the "forgetting" functors $\FF(Y,\a)\ra\D^b(Z/H)$ and
$\FF(Z,\b)\ra\D^b(Y/H)$ described above. Let
$\ov{Q}:\FF(Z,\b)\ra\FF(Y,\a)$ be the functor
defined in the same way as $\ov{R}$
but with $Y$ and $Z$ interchanged. Then it is compatible with
the ``inverse'' Fourier transform $\D^b(Y/H)\ra\D^b(Z/H)$ given by the kernel
$\ov{L}_{y,z}\simeq\ov{L}_{z,y}^{-1}$. Its composition with
the direct Fourier transform is isomorphic to a shift in the derived
category and it is easy to see that this isomorphism extends to our
additional  data,  so  that  $\ov{Q}$  is
quasi-inverse to $\ov{R}$ up to shift.
\ed

Consider the following example. Let $X=\hat{A}\times A$,
$L=L_A$ be the standard split symplectic biextension,
$Y=\hat{A}\sub\hat{A}\times A$ be its standard lagrangian
subvariety. Let $f:A\ra \hat{A}$ be a symmetric morphism,
$Z=A/\ker(f_n)\simeq (f,n\id_A)(A)\sub\hat{A}\times A$
where $f_n=f|_{A_n}$. Then it is
easy to see that $Z$ is lagrangian. Assume in addition that
$mn\ker(f)=0$ for some $m$ relatively prime to $n$.
Then we claim that the projection $X\ra X/Z$ splits. Indeed,
changing $m$ if necessary we may assume that $m+kn=1$ for
some integer $k$. Notice that we have an isomorphism
$X/Z\simeq \hat{A}/f(A_n)\simeq A/(n\ker(f))$.
Now we can define the splitting morphism
$$X/Z\simeq A/(n\ker(f))\stackrel{(k,m)}{\ra} A/\ker(f)\times
A\simeq X.$$
Hence, $e_Y=e_Z=0$ and we get an equivalence of derived
categories
$$\D^b(A)=\D^b(X/Y)\simeq\D^b(X/Z)\simeq\D^b(Z)=
\D^b(A/ker(f_n))$$
--- here we used the Fourier-Mukai equivalence for
$Z$ and $X/Z\simeq\hat{Z}$.

The proof of Theorem \ref{main} shows that we can eliminate the assumption
that there exists a biextension $P$ of $X^2$ such that
$L\simeq P\ot\si^*P^{-1}$ once we can define the categories in question
without it.
In fact, if the characteristic of the ground field is not equal to 2,
we can do it as follows. Let $Y\sub X$ be a lagrangian subvariety.
Choose another lagrangian subvariety $Z\sub X$ such that $Y\cap Z$
is finite. Then we have a $\G_m$-valued symplectic form on
$Y\cap Z$ defined by the canonical duality between
$Y\cap Z=\ker(Y\ra X/Z)=\ker(Y\ra\hat{Z})$ and
$\ker(Z\ra\hat{Y})=\ker(Z\ra X/Y)=Y\cap Z$. Since the characteristic is
different from 2 there exists a central extension $G$ of $Y\cap Z$
by $\G_m$  with  such  commutator  form  (unique  up  to  an
isomorphism), so we have the corresponding class
$e_Y\in H^2(X/Y,\G_m)$ which doesn't depend on the choices made
(and coincides with the one defined previously using $P$).

\begin{thm} Assume that the characteristic of the ground field is different
from 2. Then for every pair of lagrangian subvarieties $Y$ and $Z$
the categories
$\D^b(X/Y,e_Y)$ and $\D^b(X/Z, e_Z)$ are equivalent.
\end{thm}

\Pf . As before we may assume that $Y\cap Z$ is finite.
Now choose a central extension $G$ of $Y\cap Z$ inducing the canonical
symplectic form on it and define the functor
from the category of weight-1 $G$-equivariant objects of $\D^b(Z)$ to that
of weight-1 $G^{-1}$-equivariant objects of $\D^b(Y)$ by the formula
$$\ov{R}(\FF)_y=\int_{Z/H} L_{z,y} \FF_z dz$$
where $H\sub Y\cap Z$ is lagrangian. As above it is easy to check that
this is an equivalence.
\ed

\begin{rem}
The constructed equivalences are not canonical and they don't
agree for triples of lagrangian subvarieties. The corresponding
analogue of Maslov index and a partial generalization of this
theory to abelian schemes will be discussed in a forthcoming
paper.
\end{rem}

\bigskip

\centerline{\bf Appendix. Descent for derived categories}

\bigskip

An unpleasant property of the derived category of
coherent sheaves of $\O_S$-modules on a scheme $S$ is that
one can not glue this category from
its counterparts over open parts of $S$.
However, the following descent result holds.

\vspace{3mm}

\noindent
{\bf Theorem A.} {\it
Let $p:S'\ra S$ be a finite flat morphism.
Then the category $\D^b(S)$ is equivalent to the
the following category $\D^b(S',p)$:
its objects are pairs $(\FF,f)$ where
$\FF\in\D^b(S')$, $f:p_1^*\FF\wt{\ra} p_2^*\FF$
is an isomorphism in $\D^b(S'\times_S S')$
(where \break  $p_i:S'\times_S S'\ra   S'$,   $i=1,2$   are  the
projections) satisfying the following cocycle condition
$p_{23}^*f\circ   p_{12}^*f=p_{13}^*f$   over
$S'\times_S S'\times_S S'$.}

\vspace{3mm}

\Pf . Let $p^*:\D^b(S)\ra\D^b(S',p)$ be the natural functor.
Let us check first that $p^*$ is fully faithful.
Assume that we have a morphism $f:p^*\FF\ra p^*\GG$ in $\D^b(S')$
such that the following diagram is commutative:

\begin{equation}
\setlength{\unitlength}{0.20mm}
\begin{array}{ccccc}
p_1^*p^*\FF & \lrar{} & p_2^*p^*\FF  \\
\ldar{p_1^*(f)} & & \ldar{p_2^*(f)} \\
p_1^*p^*\GG & \lrar{} & p_2^*p^*\GG
\end{array}
\end{equation}

Applying the functor $p_{1*}$ to this diagram and composing it
with the adjunction morphism $p^*\FF\ra p_{1*}p_1^*p^*\FF$
we get the following diagram

\begin{equation}
\setlength{\unitlength}{0.20mm}
\begin{array}{ccccc}
p^*\FF & & \\
\ldar{f} & \ldrar{} & \\
p^*\GG &\lrar{} & p^*p_*p^*\GG
\end{array}
\end{equation}
where the diagonal morphism is $p^*f'$, $f':\FF\ra p_*p^*\GG$
is obtained from $f$ by adjunction.
Let us denote by $f''$ the composition
$$\FF\stackrel{f'}{\ra}p_*p^*\GG\simeq p_*\O_{S'}\ot\GG\ra
(p_*\O_{S'}/\O_S)\ot\GG.$$
Then it follows from the diagram above that $p^*(f'')=0$.
Since $p_*p^*(?)\simeq p_*\O_{S'}\ot ?$ it follows that
$f''=0$, hence $f'$ factors through a morphism
$\ov{f}:\FF\ra\GG$ and $f=p^*(\ov{f})$. Thus, the functor
$p^*:\D^b(S)\ra\D^b(S',p)$ is full and faithful. It
remains to check that any object of $\D^b(S',p)$ belongs to
its essential image. This is easy to prove by devissage with
respect to the standard $t$-structure on $\D^b(S')$ since
the corresponding truncation functors are compatible
with descent data (the base of induction is provided by
the classical descent for coherent sheaves).
\ed


\begin{thebibliography}{99}
\bibitem{Breen} L. Breen, {\it Functions theta et th\'eor\`eme du cube},
Lecture Notes in Math., vol. 980, Springer-Verlag, 1983.
\bibitem{Des} P. Deligne, {\it Le symbole mod\'er\'e}, Publ. Math. IHES 73
(1991), 147--181.
\bibitem{Mum} D. Mumford, {\it Abelian varieties}.
Second edition. Oxford University Press. 1974.
\bibitem{Weilrep} A. Polishchuk, {\it Biextensions,
Weil representation on derived categories and theta-functions.} Preprint.
\end{thebibliography}
\end{document}